\def\datum{May 13, 2009}
\documentclass[12pt]{amsart}
\usepackage{a4,amscd,amsmath,amssymb,amsfonts}
\usepackage[all]{xy}
\usepackage{cancel}
\def\sO{\mathcal O}
\def\sD{\mathcal D}

\def\sF{\mathcal F}
\def\sC{\mathcal C}
\def\et{\text{\rm \'et}}
\def\rd{\text{\rm red}}

\newtheorem{thm}{Theorem}[section]
\newtheorem{prop}[thm]{Proposition}
\newtheorem{lem}[thm]{Lemma}
\newtheorem{cor}[thm]{Corollary}

\title{Two small remarks on Nori fundamental group scheme}
\author{H\'el\`ene Esnault and Ph\`ung H\^o Hai}
\address[HE]{Department of Mathematics, Univ.  Duisburg-Essen,
45127 Essen, Germany}
\email{esnault@uni-due.de}
\address[PHH]{Institute of Mathematics, P.O. Box 731, Hanoi, Vietnam}
\email{phung@math.ac.vn}
\thanks{HE is supported by the Leibniz programme of the DFG, the SFB/TR 45 of the DFG, and the ERC-Advanced Grant 226257  }
\thanks{PHH is supported in part by the Heisenberg-Grant PH-155 and the SFB/TR 45 of the DFG, and by NAFOSTED Vietnam. He also thanks Prof. J.-H. Keum for the financial support to attend the AGEA conference.}

\date{\datum}

\begin{document}
\begin{abstract} For $X$ a complete, reduced, geometrically connected scheme over a perfect
field of characteristic $p>0$, we analyze the decomposition of Nori's fundamental group scheme into its local and \'etale parts and raise the question of the relation between the geometry and the splitting of the group scheme. We also describe in categorial terms the functor which corresponds to the inclusion of the maximal reduced subgroup scheme. 
\end{abstract}
\maketitle

\section{Nori's fundamental group scheme}   Let 
$X/k$ be a complete, reduced, geometrically connected scheme over a perfect
field  $k$. 
Let us briefly recall Nori's construction\cite{Nori1}  of the fundamental group scheme of $X$.
 A vector bundle $V$ on $X$ is said to be finite if it satisfies
a non-trivial polynomial equation with integral coefficients $f(V)\cong g(V)$ for two polynomials $f, g\in {\mathbb N}[X], f\neq g$. Here  $m\cdot V=V\oplus \ldots \oplus V$ ($m$-times) and $X^m(V)=V^{\otimes m}$. 
Subquotients of finite bundles are called essentially finite bundles. Nori showed that when $f: C\to X$ is a morphism of a smooth projective curve to $X$, then $f^*V$  is semi-stable of degree zero in the sense of Mumford whenever $V$ is an essentially finite bundle. Essentially finite bundles form an abelian rigid tensor category  $k$-linear category $\sC^N(X)$, where  morphisms are morphisms of 
vector bundles. We call these bundles Nori finite
bundles. For $X$ fixed, we shorten the notation by $\sC^N:=\sC^N(X)$.

Let us fix $x\in X(k)$. The fiber functor $V\mapsto V_{|x}$ endows 
$\sC^N$ with the structure of 
a neutral Tannaka category. Tannaka duality yields an affine (profinite) group
scheme $\pi^N(X,x)$, called {\it Nori fundamental group scheme of} $X$ {\it with base point}  $x$. For $x\to X$ fixed, we shorten the notation by $\pi^N:=\pi^N(X,x)$. 

For each $V\in \sC^N$, we denote by $\langle V\rangle$ the full subcategory of
all subquotients of direct sums of tensor powers of $V$. This is a  full
Tannaka subcategory
of $\sC^N(X)$ (with the same fiber functor). The Tannaka group of this category
is denoted by $G(V,x)$, or $G(V)$ for short. It is a finite group scheme. There is a canonical surjection 
$\pi_V:\pi^N\to G(V)$ and $\pi^N$ is  the projective limit of $G(V)$.
Furthermore, Tannaka duality applied to $\langle V\rangle$ also yields a $G(V)$-principal
bundle $p_V:Y_V\to X$ with the property that $Y_V$ is connected, is endowed with a rational point above $x$ and $p_V{}^*(W)$ is trivial for any
$W\in \langle V\rangle$.

Let
$\sC^\et$ (resp. $\sC^F$) be  the subcategory of bundles $V$ in $\sC^N$, such that $G(V)$ is
\'etale (resp. local). Tannaka duality applied to $\sC^\et$ (resp. $\sC^F$) and the fiber functor at $x$ 
yields the group scheme  $\pi^\text{\rm \'et}(X,x)$ (resp. $\pi^F(X,x)$). For $x\to X$ fixed,  we 
shorten the notation by $\pi^\et:= \pi^\text{\rm \'et}(X,x)$ (resp. $\pi^F:=\pi^F(X,x)$).
The inclusion functors $C^\text{\rm \'et}\subset \sC^N$ and $\sC^F\subset \sC^N$
yield surjective flat homorphisms  $\pi^N \to \pi^\text{\rm \'et} $
and $\pi^N \to \pi^F $. 

The group scheme $\pi^\et$ is pro-\'etale, while the group scheme $\pi^F$ is pro-local.
In fact $\pi^\et$ (resp. $\pi^F$) is the largest pro-\'etale (resp. pro-local) quotient of $\pi^N$.
In~\cite{EHX} the relationship between these    group schemes has been studied.
It is shown in particular that there is a canonical homorphism
$\pi^N\to \pi^\et\times_k \pi^F$ which is flat surjective but generally not an isomorphism.
The description of the kernel of this map in terms of Tannaka duality was given.

For a finite group scheme $G$ over $k$, let $G^0$ denote the connected component
of the unit element. The reduced subscheme  $G_\rd \subset G$ is a 
subgroup scheme. The composite homomorphism $G_\rd\to G\to G^\et$ is an isomorphism and furthermore,  $G^0$ is a normal
subgroup. So  $G_\rd\cong G/G^0$. In other words, $G$ is the semi-direct product of $G^0$ with $G_\rd$, with $G^0\subset G$ normal. The aim
of this short note is to discuss the behavior of this semi-direct presentation of $G$ in the pro-system defining   $\pi^N$.
\section{The reduced part of $\pi^N$}
\subsection{The pro-local group $(\pi^N)^0$} Let $(\pi^N)^0$ be the pro-local subgroup scheme of $\pi^N$. That is, with the notations as in the previous section
\begin{equation} (\pi^N)^0=\varprojlim_{V\in\sC^N(X)}G(V)^0.\end{equation}

\begin{lem} $(\pi^N)^0$ is the kernel
of $\pi^N\to \pi^\text{\rm \'et}$.\end{lem}
\begin{proof} For $V \in \sC^N$ given,  $G(V)^\et= G(V)/G(V)^0$ is the largest \'etale quotient of $G(V)$. Thus by Tannaka duality, 
$G(V)^\et$ is the Tannaka group of the sub category $\langle V\rangle\cap \sC^\et \subset \langle V\rangle$.
Hence 
\begin{equation} \pi^\et=\varprojlim_{V\in\sC^N}G(V)^\et
\end{equation}
The lemma follows. \end{proof}

A description of the kernel of $\pi^N\to \pi^\et$ in terms of Tannaka duality was given
in \cite{EHX} (there it is denoted by $L$). 
It is shown that $(\pi^N)^0$ may differ from $\pi^F$, more
precisely, the composition homomorphism $(\pi^N)^0\to \pi^N \to  \pi^F$
may be not injective.

\subsection{The reduced group $\pi^N_\rd$}
In the language of function algebras, $\sO(\pi^\text{\rm \'et})$
is the largest pro-\'etale sub Hopf algebra of $\sO(\pi^N)$
and $\sO(\pi^F)$ is the largest nilpotent sub Hopf algebra
of $\sO(\pi^N)$. Furthermore, $\sO(\pi^N)$ is the ind limit of its finite dimensional sub Hopf algebras $\sO(G)$ over $k$.

On the other hand, let $\frak{N}$ be the nilradical of $\sO(\pi^N)$.
Then for each sub Hopf algebra $\sO(G)\subset \sO(\pi^N)$, 
$\frak{N}\cap \sO(G)$ is the nilradical of $\sO(G)$. 
Further we have
$\sO(G)/(\frak{N}\cap \sO(G))=\sO(G_\text{\rm red})$, where $G_\text{\rm red}$ is the
reduced subgroup of $G$ with the same underlying topological space. Thus the
quotient $\sO(\pi^N)/\frak{N}$ is the ind-limit of  $\sO(G_\text{\rm red})$,
its spectrum $\pi^N_\text{\rm red}$ is the largest pro-\'etale 
subgroup of $\pi^N$, which is also the pro-reduced subscheme of $\pi^N$.

 \begin{prop}With the settings above we have
 \begin{itemize}\item[(i)] The composition homomorphism
 \begin{equation}\pi^N_\text{\rm red}\to \pi^N\to \pi^\text{\rm \'et}\end{equation}
is an isomorphism. Consequently the inverse of this map
defines a splitting of $\pi^N\to \pi^\text{\rm \'et}$.
\item[(ii)] $\pi^N$ is  the semi-direct product of its subgroups
$(\pi^N)^0$ and $\pi_\text{\rm red}$, with $(\pi^N)^0$ normal.
  \end{itemize}
 \end{prop}

\begin{cor} The composition map
$(\pi^N)^0\to\pi^N\to \pi^F$ is an isomorphism if and only if $\pi^N$ is the
direct product of $\pi^F$ with $\pi^\text{\rm \'et}$ in a compatible way in the  pro-system, i.e.  
$$\pi^N=\varprojlim_V G(V)^F\times_k G(V)^\et$$
\end{cor}
\begin{proof}
For a finite group scheme $G$, the claim of Corollary holds. Passing to limit in the
prosystem defining $\pi^N$ we obtain the claim of Corollary.
\end{proof}


If $G$ is a commutative (pro)-finite group scheme, then $G=G^0\times_k G^\et$.
 Of course $G$ need not
be commutative to split in this way. We raise the question of the geometric conditions on $X$ which force an isomorphism $\pi^\et\cong \pi^F\times_k \pi^N_\text{\rm red}$. More precisely we ask for the relation between this strong splitting condition and the commutativity of $\pi^N$. If $X$ has dimension 1,  there is a simple answer.

\begin{prop} Let $X$ be a smooth projective geometrically connected curve over a char. $p>0$ algebraically closed field $k$. Let $x\in X(k)$. Then $\pi^N(X,x) \cong \pi^F(X,x)\times_k \pi^N_{\rm red}(X,x)$    if and only if $\pi^N(X,x)$ is commutative.\end{prop}

\begin{proof} 
Assume that $\pi^N=(\pi^N)^0\times_k \pi^N_\rd$. Then $(\pi^N)^0\to\pi^N\to\pi^F$
is an isomorphism. According to \cite[Thm~3.5]{EHX}, the representation category
of $(\pi^N)^0$ is equivalent to the category $\sD$, which consists of pairs $(X_S,V)$ 
where $X_S\to X$ is the principal bundle associated to a full subcategory of $\sC^\et$
generated by some object of $\sC^\et$,  morphisms are appropriately determined,
(see \cite[Defn~3.3]{EHX}), and $V\in\sC^F(X_S)$. The morphism
$(\pi^N)^0\to \pi^F$ is Tannaka dual to the functor $V\mapsto (X,V)$ (i.e. $S$ is the trivial
subcategory of $\sC^\et$). The isomorphism $(\pi^N)^0\to \pi^F$ implies that
each $V\in \sC^F(X_S)$ is the pull back of some $W\in\sC^F(X)$.

If $X$ is a smooth curve of genus $\ge 2$, Raynaud \cite{raynaud} shows that there exists an \'etale
 cyclic cover $X_L\to X$ and a $p$-torsion line bundle
on $X_L$ which does not come from a line bundle on $X$ . Thus, in this case
one cannot have isomorphism $(\pi^N)^0\to \pi^F$. On the other hand,
if $X$ has genus 1, the point $x$ gives $X$ the structure of an abelian variety. It is shown by Nori \cite{Nori2}  that $\sC^N(X)$ is 
commutative.
 Finally, if $X$ has genus $0$, then $\pi^\et=\pi^N_\text{red}=\{1\}$ so there is nothing to show. This finishes the proof. 
\end{proof}

\section{Tannaka description of the map $\pi^\text{\rm \'et}\to \pi^N$.}
Our aim in this section is to describe the functor $\text{\rm\'Et}:\sC^N\to \sC^\text{\rm \'et}$ that corresponds, through
the Tannaka duality, to the injection $\pi^\text{\rm \'et}\to \pi^N$.

\subsection{The Frobenius functor on a representation category} Let $G$ be 
a group scheme over $k$. The Frobenius functor $\sF$ on $\text{Rep}(G)$ is defined as follows.
For each representation $V$ of $G$, $\sF(V)$, as a $k$ vector space, is $V^{(1)}:=V\otimes_{F_k}k$,
where $F_k$ is the Frobenius of $k$. Let $\{e_i\}$ be a $k$ basis of $V$. If the action of $g\in G$
on $V$ is given by a matrix $(g_{ij})$, its  action  on $V^{(1)}$ is defined by the 
the matrix $(g_{ij}{}^p)$. In the dual language of functions algebras, if the coaction of
$\sO(G)$ on $V$ is $\delta:v\mapsto \sum_iv_i\otimes a_i$ then the coaction of
$\sO(G)$ on $V^{(1)}$ is 
\begin{equation}\label{coaction}\delta^{(p)}:v\otimes \lambda\mapsto \sum_i(v_i\otimes\lambda)\otimes a_i{}^p.\end{equation}

Recall that the absolute Frobenius $F: X\to X$ of $X$ induces a functor $F^*: \sC^N(X)\to \sC^N(X), \ V \mapsto F^*(V)$.
It should be a well-known fact that $\sF$ is compatible
with the Frobenius functor through the fiber functor $\omega_x$. We provide here a simple proof of this fact
(see also \cite[Thm.11]{santos}).

\begin{lem} The functor $F^*$ is compatible with the Frobenius functor on $\text{\rm Rep}(\sC^N(X))$
through fiber functor as follows. The following diagram is commutative:
\begin{equation}\label{frob}
\xymatrix{\sC^N(X)\ar[r]^{F^*}\ar[d]_{\omega_x}& \sC^N(X)\ar[d]^{\omega_x}\\
\text{\rm Rep}(G)\ar[r]_{\sF^*}&\text{\rm Rep}(G)
}
\end{equation}
\end{lem}
\begin{proof}It is easy to see that $\sF$ is a $p$-linear tensor functor.
 This means  $\sF(\lambda f)=\lambda^p\sF(f)$, for all $\lambda\in k$ and all morphisms $f$. Moreover, $\sF(V)$ can be
determined by using only ``algebraic tensor constructions'' as follows. Denote by $S^n(V)$ the
$n$-th symmetric power of $V$, i.e. the largest quotient of $V^{\otimes n}$ invariant
by all symmetries on $V^{\otimes n}$, and denote by $ST^n(V)$ the
subspace of $T^n(V)$ of symmetric tensors. Both $S^n(V)$ and $ST^n(V)$ are
representations of $G$. Moreover, the image of the composition
\begin{equation}ST^p(V)\to V^{\otimes p}\to S^p(V)\end{equation}
is naturally isomorphic to $\sF(V)=V^{(1)}$ as a $G$-representation.
Indeed, the image of $ST^p(V)$ in $S^p(V)$ is 
spanned by $\{e_i^{\otimes p}\}_{i=1,\dots,n}$
where $\{e_i\}_{i=1,\dots,n}$ is a basis of $V$
and the restriction of the coaction on $S^p(V)$
on this subspace has the same form as the action
$\delta^{(p)}$ given in (\ref{coaction}). 

On the other hand $\sF(V)$ can also be defined using only ``algebraic tensor
constructions'' as $V^{(1)}$ above.
 Since $\omega$ is
exact and compatible with the tensor structures, it
satisfies the diagram in (\ref{frob}).
\end{proof}
\begin{cor}

The functor $F^*$, restricted on $\sC^\text{\rm \'et}$, is an equivalence of categories.
\end{cor}
 
\begin{proof}  For an \'etale $k$-algebra, the absolute Frobenius homomorphism
is an isomorphism. Hence the same holds true for a pro-\'etale $k$-algebra.
Looking at the coaction in (\ref{coaction}) we see
that the Frobenius functor on $\text{Rep}(\pi^\text{\'et})$ is an equivalence
of categories. By Tannaka duality, the functor $F^*$ on $\sC^\text{\'et}$
is an equivalence of categories.\end{proof}
 We
henceforth denote by $F^{*-1}$ the inverse functor and by $F^{*-n}$ its $n$-th power, which is the inverse functor to $F^{*n}:=F^* \circ \ldots \circ F^*$ ($n$-times). 

The proof of the proposition below is now obvious.

\begin{prop} Let $\langle V\rangle$ be the full tensor subcategory of $\sC^N$, generated
by an object $V$ in $\sC^N$. The restriction of the functor \'Et to $\langle V\rangle$ is 
equivalent to functor $F^{*-n}(F^{*n}(V))$ for any $n$ larger than some $n_V$ depending
on $V$.
\end{prop}
\begin{proof}
Indeed, there exists an integer $n_V$ such that $F^{*n_V}$
is \'etale. Since $F^*$, restricted in $\sC^\et$, is an equivalence of categories,  for $n\geq n_V$,
$F^{*-n}(F^{*n}(W))$ is a well-defined, $k$-linear functor
from $\langle V \rangle \subset \sC^N$ to $\sC^\et$.

\end{proof}

\end{document}